\documentclass{article}
\usepackage{amsmath}
\usepackage{amsfonts}
\newcommand{\proof}{\textit{Proof}\textrm{. }}
\newcommand{\epr}{\hfill$\diamondsuit$\medskip}
\newtheorem{theorem}{Theorem}
\newtheorem{corollary}[theorem]{Corollary}
\newtheorem{lemma}[theorem]{Lemma}
\newtheorem{proposition}[theorem]{Proposition}
\title{On the polynomial identities of the algebra $M_{11}(E)$}
\author{Plamen Koshlukov\thanks{Partially supported by grants from CNPq
  (No. 304003/2011-5 and 480139/2012-1), and  from FAPESP  (No. 2010/50347-9)}, Thiago Castilho de Mello\thanks{Supported by PhD grant from CNPq}
\\
Department of Mathematics, IMECC, UNICAMP\\
S\'ergio Buarque de Holanda 651\\
13083-859 Campinas, SP, Brazil\\
e-mail: \texttt{plamen@ime.unicamp.br}, \texttt{tcmello@ime.unicamp.br}}
\date{}
\begin{document}
\maketitle



\begin{abstract}
Verbally prime algebras are important in PI theory. They  were described by Kemer over a field $K$ of characteristic zero: 0 and $K\langle T\rangle$ (the trivial ones), $M_n(K)$, $M_n(E)$, $M_{ab}(E)$. Here $K\langle T\rangle$ is the free associative algebra of infinite rank, with free generators $T$, $E$ denotes the infinite dimensional Grassmann algebra over $K$, $M_n(K)$ and $M_n(E)$ are the $n\times n$ matrices over $K$ and over $E$, respectively. The algebras $M_{ab}(E)$ are subalgebras of $M_{a+b}(E)$, see their definition below. The generic (also called relatively free) algebras of these algebras have been studied extensively. Procesi described the generic algebra of $M_n(K)$ and lots of its properties. Models for the generic algebras of $M_n(E)$ and $M_{ab}(E)$ are also known but their structure remains quite unclear.

In this paper we study the generic algebra of $M_{11}(E)$ in two generators, over a field of characteristic 0. In an earlier paper we proved that its centre is a direct sum of the field and a nilpotent ideal (of the generic algebra), and we gave a detailed description of this centre. Those results were obtained assuming the base field infinite and of characteristic different from 2. In this paper we study the polynomial identities satisfied by this generic algebra. We exhibit a basis of its polynomial identities. It turns out that this algebra is PI equivalent to a 5-dimensional algebra of certain upper triangular matrices. The identities of the latter algebra have been studied; these were described by Gordienko. As an application of our results we describe the subvarieties of the variety of unitary algebras generated by the generic algebra in two generators of $M_{11}(E)$. Also we describe the polynomial identities in two variables of the algebra $M_{11}(E)$.
\end{abstract}

Let $K$ be a field of characteristic 0, and denote by $E$ the infinite dimensional Grassmann algebra over $K$. If $V$ is a $K$-vector space with a basis
$e_1$, $e_2$, \dots{} then $E=E(V)$ is the vector space with a basis consisting of 1 and all products $e_{i_1}e_{i_2}\cdots e_{i_k}$, $i_1<i_2<\cdots<i_k$,
$k\ge 1$. The multiplication in $E$ is induced by the anticommutative law for the $e_i$'s, that is by $e_ie_j=-e_je_i$ for all $i$ and $j$. The Grassmann
algebra plays an extremely important role in the theory of PI algebras. It is one of the natural examples of a PI algebra that satisfies no standard
identity whenever char$K=0$. (Note that every PI algebra over a field of characteristic $p>0$ satisfies some standard identity, due to a result of Kemer
\cite{kemer}). The Grassmann algebra is $\mathbb{Z}_2$-graded: one can verify that $E=E_0\oplus E_1$ where the vector subspaces $E_i$ are spanned by all
elements from its basis such that $k\equiv i\pmod{2}$. Then $E_aE_b\subseteq E_{a+b}$ where the latter sum is modulo 2. Recall that $E_0$ is just the
centre of $E$. This gives $E$ the structure of a 2-graded algebra (also called superalgebra). One of the most important developments in PI theory came with
Kemer's results on the structure of the T-ideals. These results led Kemer directly to the positive solution of the long standing Specht problem, among
others (see \cite{kemerbook}). Kemer based his theory on the description of the T-prime (also called verbally prime) algebras. These turn out to be the
matrix algebras $M_n(K)$, the algebras $M_n(E)$ of the matrices with entries in $E$, and the algebras $M_{ab}(E)$. The latter are  subalgebras of
$M_{a+b}(E)$; they consist of all block matrices with blocks of sizes $a\times a$ and $b\times b$ on the main diagonal with entries from $E_0$, and the
remaining, off-diagonal blocks with entries from $E_1$. Afterwards Kemer proved that every finitely generated PI algebra satisfies the same identities as a
suitable finite dimensional algebra. Moreover, if $A$ is any PI algebra then it satisfies the same polynomial identities as the Grassmann hull of a
suitable finite dimensional superalgebra. Recall that if $A=A_0\oplus A_1$ is 2-graded then its Grassmann hull is $G(A)=A_0\otimes E_0\oplus A_1\otimes
E_1$. The reader can find the details about Kemer's results in his monograph \cite{kemerbook}.

In spite of the importance of the T-prime algebras not much is known about the concrete form of their polynomial identities. The T-ideals of $M_n(K)$ are
known only for $n=1$ and 2, see \cite{razmal, dral}, and also \cite{razmbook} for a streamlined version of the proof. When $K$ is infinite and of
characteristic different from 2, a basis of the identities of $M_2(K)$ is also known, see \cite{kjam2, jcpk}. When $K$ is a finite field, bases of the
identities of $M_n(K)$ are known for $n\le 4$, see \cite{maltsev, genov, genovsid}, respectively. It should be noted that the methods applied in the case
of $K$ finite are completely different from the ones when $K$ is infinite. The identities of the Grassmann algebra $E$ were described in \cite{krreg} in
characteristic 0, and by various authors over any field, see the bibliography of \cite{gk}. In \cite{popov} Popov described the identities of $E\otimes E$
over a field of characteristic 0. Note that according to Kemer's theory $E\otimes E$ and $M_{11}(E)$ satisfy the same polynomial identities, that is they
are PI equivalent. Note also that these two algebras are not PI equivalent in positive characteristic $p>2$, see \cite{afk}.

When one studies the identities of $M_n(K)$ the algebra of the generic matrices comes into play. Procesi made an extensive use of the algebra of generic matrices and described very many of its properties, see for example \cite{Procesi}. Berele in \cite{bereleca} constructed generic algebras for the remaining T-prime algebras. We recall these constructions below.

In an earlier paper \cite{pktcm} we started studying the relatively free algebra of $M_{11}(E)$ in two generators. We were able to give quite precise description of its centre. We proved that the centre of this algebra is a direct sum of the field and a nilpotent ideal of the algebra, and we showed that the centre contains quite a lot of non-scalars. In this way we gave an answer to a question posed in \cite{bereleca}. Below we recall results of \cite{pktcm} needed for our exposition. Our notation is similar to that of \cite{pktcm}.

In this paper we work over a fixed field of characteristic 0. We describe the polynomial identities of the generic algebra in two generators for $M_{11}(E)$. It turns out the basis of its identities consists of three polynomials given in explicit form. It is worth mentioning that the same polynomials form a basis of the identities of a certain subalgebra of the $3\times 3$ upper triangular matrices. This algebra appears in the study of extremal properties of the codimension sequence, see for example \cite{smav}. Later on its identities were described by Gordienko, see for example \cite{gordienko}. Our methods are based on the representation theory of the symmetric and the general linear groups, on the results from \cite{pktcm}, and on the description of the polynomial identities of $M_{11}(E)$ given by Popov in \cite{popov}. Moreover we describe the unitary subvarieties of the variety generated by our algebra. We also prove that asymptotically the identities of every proper (unitary) subvariety behave exactly as the identities satisfied by the algebra of upper triangular matrices of order two, or as those of commutativity.

As an application of the techniques developed here we describe also the polynomial identities of $M_{11}(E)$ in two variables. To this end we use  some ideas from a paper by Nikolaev, \cite{nikolaev} where the identities in two variables of $M_2(K)$ were described.

This paper is an extended version of another paper submitted for publication.

\section{Preliminaries}
Throughout we consider unitary associative algebras over a field $K$ of characteristic 0. Let $K\langle T\rangle$ be the free associative algebra freely generated over $K$ by the set $T=\{t_1,t_2, \dots\}$. If $A$ is a PI algebra we denote by $I=T(A)\subseteq K\langle T\rangle$ its T-ideal, that is the ideal of all identities of $A$. The algebra $U(A)=K\langle T\rangle/I$ is the relatively free (or generic) algebra in the variety of algebras defined by $A$. When $T$ is a finite set, say $T=\{t_1,\dots,t_k\}$ one obtains the relatively free (or generic) algebra of $A$ of rank $k$ and denotes it by $U_k(A)$. We shall use the same letters $t_i$ for the generators of $K\langle T\rangle$ and for their images under the canonical projection $K\langle T\rangle\to K\langle T\rangle/T(A) = U(A)$.

We recall the construction of the free supercommutative algebra $K[X;Y]$. Let $K\langle X\cup Y\rangle$ be the free associative algebra freely generated by the set $X\cup Y$ where $X\cap Y=\emptyset$. This algebra is 2-graded in a natural way assuming the variables in $X$ of degree 0, and those in $Y$ of degree 1. Let $I$ be the ideal generated by all $ab-(-1)^{|a|\cdot|b|}ba$ where $a$ and $b$ run over the homogeneous elements in $K\langle X\cup Y\rangle$, and $|a|$ is the $\mathbb{Z}_2$-degree of the homogeneous element $a$, and put $K[X;Y]= K\langle X\cup Y\rangle/I$. Then $K[X;Y]\cong K[X]\otimes_K E(Y)$, here $E(Y)$ is the Grassmann algebra on the vector space with a basis $Y$, see for more detail \cite{bereleca}. If a 2-graded algebra $A=A_0\oplus A_1$ satisfies $ab-(-1)^{|a|\cdot|b|}ba=0$ for all homogeneous $a$ and $b$ then it is called supercommutative. Clearly the Grassmann algebra is supercommutative.

Take $X=\{x_{ij}^r\}$, $Y=\{y_{ij}^r\}$ where $1\le i,j\le n$, $r=1$, 2, \dots; here $r$ is an upper index, not an exponent. One defines
the matrices $A_r=(x_{ij}^r)$, $B_r=(x_{ij}^r+y_{ij}^r)$, $C_r=(z_{ij}^r)$ where $z=x$ whenever $1\le i,j\le a$ or $a+1\le i,j\le a+b$, and $z=y$ for all remaining possibilities for $i$ and $j$. Suppose $a+b=n$, and consider the following subalgebras of $M_n(K[X;Y])$. The first is generated by the generic matrices $A_r$, $K[A_r\mid r\ge 1]$. It is isomorphic to the relatively free (or universal) algebra $U(M_n(K))$ of $M_n(K)$. In
\cite[Theorem 2]{bereleca} it was proved  that $U(M_n(E))\cong K[B_r\mid r\ge 1]$, also  $U(M_{ab}(E))\cong K[C_r\mid r\ge 1]$. The relatively
free algebras of finite rank $k$, denoted by $U_k$, can be obtained by letting $r=1$, \dots, $k$, that is by taking the first $k$ matrices.

The algebra $K\langle T\rangle$ is multigraded, counting the degree of its monomials in each variable.
We work over a field $K$ of characteristic 0 therefore every T-ideal is generated by its multilinear elements, see for example \cite[Section 4.2]{drenskybook}.

The polynomial identities of $M_{11}(E)$ were described by Popov in characteristic 0, see the main theorem of \cite{popov}. The theorem reads that the polynomials
\begin{equation}
\label{basispim11}
[[t_1,t_2]^2,t_1], \qquad [[t_1,t_2],[t_3,t_4],t_5]
\end{equation}
generate the T-ideal of $E\otimes E$, and of $M_{11}(E)$ as well.
Here $[a,b]=ab-ba$ is the usual commutator of $a$ and $b$. The higher commutators without inner brackets will be considered left normed that is $[a,b,c] = [[a,b],c]$, and so on.

Let $L(T)$ be the free Lie algebra on the free generators $T$. If one substitutes the usual product in an associative algebra $A$ by the bracket $[a,b] = ab-ba$ one gets a Lie algebra denoted by $A^-$. It is well known that $L(T)$ is the subalgebra of $K\langle T\rangle^-$ generated by $T$. Moreover $K\langle T\rangle$ is the universal enveloping algebra of $L(T)$. Choose an ordered basis of $L(T)$ consisting of $T$ and left normed commutators. Suppose further that if $u$ and $v$ are elements of the basis then $\deg u<\deg v$ implies $u<v$, in this way the free generators in $T$ precede all remaining basis elements. Then a basis of $K\langle T\rangle$ is given by 1 and all elements $t_1^{n_1}\cdots t_k^{n_k} u_1\cdots u_m$ where $n_i\ge 0$, and $u_i$ are commutators, $u_1\le\cdots \le u_m$. Let $B(T)$ be the subalgebra of $K\langle T\rangle$ generated by 1 and all commutators of degree at least two. Thus $B(T)$ is spanned by 1 and the products of commutators. The elements of $B(T)$ are called proper polynomials. As we work with unitary algebras it is well known that every T-ideal $I$ is generated by its proper elements, see for example \cite[Section 4.3]{drenskybook}.

\section{The generic algebra of $M_{11}(E)$ in two generators}
In the paper \cite{pktcm} we studied the generic algebra of $M_{11}(E)$ in two generators, $F=K[C_1,C_2]$  where $C_1=\begin{pmatrix}
x_1&y_1\\ y_1'& x_1'\end{pmatrix}$, $C_2=\begin{pmatrix} x_2&y_2\\ y_2'& x_2'\end{pmatrix}$. The entries of $C_1$ and $C_2$ lie in the free supercommutative algebra $K[X;Y]$, $X=\{x_1,x_2,x_1', x_2'\}$, $Y=\{y_1,y_2, y_1', y_2'\}$. Here we recall some results of \cite{pktcm} that we need.

The algebra $K[X;Y]$ is 2-graded in a natural way. It can be given a $\mathbb{Z}$-grading by counting the degree in the variables of $Y$ only: $K[X;Y] = \oplus_{n\in\mathbb{Z}}
K[X;Y]^{(n)}$ where $K[X;Y]^{(n)}$ is the span of all monomials of degree $n$ in the variables from $Y$. Obviously $K[X;Y]^{(n)}=0$ unless $0\le n\le 4$. Also
\begin{eqnarray*}
K[X;Y]_0 &=& K[X;Y]^{(0)} + K[X;Y]^{(2)} + K[X;Y]^{(4)}; \\
K[X;Y]_1 &=& K[X;Y]^{(1)} + K[X;Y]^{(3)}.
\end{eqnarray*}
We set $B_0=\{1\}$, $B_1=\{y_1,y_2,y_1', y_2'\}$, $B_2=\{y_1y_2, y_1y_1', y_1y_2', y_2y_1', y_2y_2', y_1'y_2'\}$,
$B_3=\{y_1y_2y_1', y_1y_2y_2', y_1y_1'y_2', y_2y_1'y_2'\}$,  $B_4=\{y_1y_2y_1'y_2'\}$. Then for $n=0$, 1, 2, 3, 4, $K[X;Y]^{(n)}$ is a free module over $K[X]$, with a basis $B_n$. As a consequence $K[X;Y]$ is a free module over $K[X]$ with a basis $B=B_0\cup B_1\cup B_2\cup B_3\cup B_4$. Also the ideals in $K[X;Y]$ are $K[X]$-submodules.

The following elements were introduced in \cite{pktcm}.
\begin{eqnarray*}
h_1&=& y_1y_2y_1'y_2';\\
h_2 &=& y_1y_2 (y_1'(x_2'-x_2) - y_2'(x_1'-x_1)); \\
h_3 &=& y_1'y_2' (y_1(x_2'-x_2) - y_2(x_1'-x_1)); \\
h_4 &=& (y_1'(x_2'-x_2) - y_2'(x_1'-x_1))
(y_1(x_2'-x_2) - y_2(x_1'-x_1)).
\end{eqnarray*}
It is immediate to check they satisfy the following relations in $K[X;Y]$.
\begin{eqnarray*}
&&h_1y_1 = h_1y_1' = h_1y_2 = h_1y_2' = 0; \quad h_2y_1 = h_2y_2 = h_3y_1' = h_3y_2' = 0;\\
&&h_2y_1' = h_3y_1 = (x_1'-x_1)h_1; \quad h_2y_2' = h_3y_2 = (x_2'-x_2) h_1; \\
&&h_4y_1 = (x_1'-x_1)h_2; \quad h_4y_2 = (x_2'-x_2) h_2; \\
&&h_4y_1' = -(x_1'-x_1) h_3; \quad h_4y_2' = -(x_2'-x_2)h_3.
\end{eqnarray*}
We shall also need the elements
\begin{eqnarray*}
&&q_n(x_1,x_1') = \sum\nolimits_{i=0}^nx_1^i x_1'^{n-i}; \quad
Q_n(x_2,x_2') = q_n(x_2,x_2');\\
&&r_n(x_1, x_1') = \sum\nolimits_{i=0}^{n-1} (n-i)x_1^{n-1-i}x_1'^i; \quad
R_n(x_2,x_2') = r_n(x_2,x_2');\\
&&s_n(x_1,x_1') = r_n(x_1', x_1); \quad
S_n(x_2,x_2') =  s_n(x_2,x_2').
\end{eqnarray*}
One verifies by an obvious induction that
\begin{eqnarray*}
&&r_n=q_{n-1}+x_1r_{n-1}; \quad s_n=q_{n-1} + x_1's_{n-1}; \quad s_n+r_n = (n+1) q_{n-1}; \\
&&q_n=x_1^n + x_1' q_{n-1} = x_1'{}^n + x_1q_{n-1};\quad (x_1'-x_1)q_{n-1} = x_1'{}^n - x_1^n;\\
&& x_1^nx_1'{}^m - x_1^mx_1'{}^n = (x_1'-x_1)(q_nq_{m-1} - q_mq_{n-1}).
\end{eqnarray*}
We compute directly that for every $m$ and $n$ we have
\begin{eqnarray*}
C_1^n &=&\begin{pmatrix}
x_1^n+y_1y_1'r_{n-1}& y_1q_{n-1}\\
y_1'q_{n-1}& x_1'{}^n-y_1y_1's_{n-1}
\end{pmatrix}; \\
C_2^m &=& \begin{pmatrix}
x_2^m + y_2 y_2' R_{m-1} & y_2 Q_{m-1}\\
y_2' Q_{m-1} & x_2'{}^m - y_2y_2' S_{m-1}
\end{pmatrix}.
\end{eqnarray*}
In this way the product $C_1^nC_2^m$ equals
\begin{equation}
\label{explicitproduct}
C_1^nC_2^m = \begin{pmatrix}
x_1^nx_2^m+a+d & y_1x_2'{}^mq_{n-1} + y_2x_1^n Q_{m-1} + c\\
y_1'x_2^m q_{n-1} + y_2'x_1'{}^nQ_{m-1} + c' & x_1'{}^n x_2'{}^m + a'+d'
\end{pmatrix}
\end{equation}
where $a$, $a'\in K[X;Y]^{(2)}$, $d$, $d'\in K[X;Y]^{(4)}$, and $c$, $c'\in K[X;Y]^{(3)}$.

Let $A=\begin{pmatrix} a&b\\ c&d\end{pmatrix}\in F\cong K[C_1,C_2]$ be central. It was shown in \cite{pktcm} that \[
d-a\in J= Ann ((x_2'-x_2)y_1 - (x_1'-x_1)y_2)\cap
Ann ((x_2'-x_2)y_1' - (x_1'-x_1)y_2').
\]
By \cite[Proposition 5]{pktcm} the $K[X]$-module $J$ is spanned by $\{h_1,h_2,h_3,h_4\}$. Corollary 6 from \cite{pktcm} states that the matrix $A$ commutes with $C_1$ and $C_2$ if and only if $b=f_4 h_2$, $c = -f_4 h_3$, $d = a+f_1h_1+f_4h_4$ for some $f_1$, $f_4\in K[X]$. Hence $A$ is central if and only if
$
A= aI + f_1\begin{pmatrix} 0&0\\ 0&h_1\end{pmatrix} + f_4\begin{pmatrix} 0&h_2\\ -h_3& h_4\end{pmatrix}$.

Following the notation used in \cite{pktcm} we define matrices in $F$:
\[
A_0=\begin{pmatrix} h_1&0\\ 0&h_1\end{pmatrix}; \quad
A_1=\begin{pmatrix} 0&0\\ 0&h_1\end{pmatrix}; \quad
A_2=\begin{pmatrix} 0&h_2\\ -h_3&h_4\end{pmatrix}; \quad
A_3=\begin{pmatrix} h_4&0\\ 0&h_4\end{pmatrix}.
\]
An element $a\in F$ is \textsl{strongly central} if it is central, and also for each $b\in F$ the element $ab$ is central in $F$.
One checks (\cite[Lemma 7]{pktcm}) that for arbitrary $\alpha_i\in K[X]$ the elements $\alpha_0A_0 + \alpha_1A_1 + \alpha_2A_2 + \alpha_3A_3$ are strongly central in $F$.

We shall need a couple of technical statements from \cite{pktcm}. Take a left normed commutator $f(t_1,t_2) = [t_1,t_2,t_{i_3}, \ldots, t_{i_k}]$, $i_j=1$, 2, such that $\deg_{t_1} f=n$, $\deg_{t_2} f=m$, $n+m=k$. Then Lemma 8 from \cite{pktcm} states that $f(C_1,C_2) = (x_1'-x_1)^{n-1} (x_2'-x_2)^{m-1} A(k)$ where
\begin{eqnarray*}
A(k) &=& \begin{pmatrix} F(k) & y_1(x_2'-x_2) - y_2 (x_1'-x_1)\\
(-1)^k (y_2'(x_1'-x_1) - y_1'(x_2'-x_2))& F(k) \end{pmatrix},\\
F(k) &=& \frac{(y_1(x_2'-x_2) - y_2(x_1'-x_1)) y_i' + (-1)^k y_i (y_2'(x_1'-x_1) - y_1'(x_2'-x_2))}{x_i'-x_i} .
\end{eqnarray*}
In the expression for $F(k)$ the index $i$ stands for $i_k$. The  formula for $f(C_1,C_2)$ yields that if $f_1$, $f_2$ are  non-zero commutators in $F$
then $f_1(C_1,C_2)f_2(C_1,C_2)$ is strongly central in $F$, see \cite[Lemma 9]{pktcm}.

\section{Some identities for the algebra $F$}

At the beginning of the previous section we denotted by $F$ the $K$-algebra generated by the matrices $C_1$ and $C_2$, $F=K[C_1, C_2]$.

\begin{lemma}
\label{3comm}
Let $f_1(t_1,t_2)$, $f_2(t_1,t_2)$, $f_3(t_1,t_2)$ be three commutators and put $g=f_1f_2f_3$. Then $g(C_1,C_2)$ vanishes on $F$.
\end{lemma}

\proof
Write the $f_i$ as linear combinations of left-normed commutators. It suffices to consider the case when all $f_i$ are left-normed. By \cite[Lemma 9]{pktcm} the product $f_1(C_1,C_2) f_2(C_1,C_2)$ is, in turn, a combination of the matrices $A_0$, $A_2$, $A_3$ defined above. Now the entries of the matrices $A_i$ are either zeros or, up to a sign, some of the elements $h_j$, $1\le j\le 4$. Looking at the above expression for $f_3(C_1,C_2)$ we see that its entries vanish when multiplied by $h_j$.
\epr

\begin{corollary}
\label{iffproper}
Let $f(t_1,t_2) = t_1^nt_2^m u_1^{k_1}\cdots u_r^{k_r}$ where the $u_i$ are left-normed commutators in $t_1$ and $t_2$. Then $f(C_1,C_2)=0$ in $F$ if and only if $k_1+\cdots+k_r\ge 3$.
\end{corollary}

\proof Lemma~\ref{3comm} implies the ``if" part. Suppose  $k_1+\cdots+k_r\le 2$. According to the proof of \cite[Lemma 9]{pktcm} the product of two
left-normed commutators, evaluated on $C_1$ and $C_2$ in $F$, is a linear combination of the matrices $A_0$, $A_2$ and $A_3$, the latter with coefficient $\pm 1$, 
hence it cannot be 0 in $F$. Also the matrices $C_1$ and $C_2$ are not zero divisors in $F$ according to \cite[Lemma 2]{pktcm}. \epr

We make use of the following two equalities that are valid in every associative algebra. They are well known and their proofs consist of an easy induction. We separate them into a lemma for further reference.

\begin{lemma}
\label{equalities}
Let $z$, $a$, $b\in K\langle T\rangle$, and $n\ge 1$ then
\[
[z^n,b] = \sum_{i=0}^{n-1} z^{n-i-1} [z,b] z^i; \qquad
[a,b]z^n = \sum_{i=0}^n {n\choose i} z^i [a,b,\underbrace{z, \ldots,z}_{n-i}].
\]
\end{lemma}

\begin{lemma}
\label{leftright}
Let $u$ be a left-normed commutator in $K\langle T\rangle$ and let $v\in K\langle T\rangle$. Then $uv = \sum_i v_i u_i$ where $u_i$ are left-normed commutators, $v_i\in K\langle T\rangle$, and $\deg u_i\ge \deg u$ for all $i$.
\end{lemma}

\proof
It suffices to consider $v=t_1$, then $ut_1 = [u,t_1] + t_1u$. Then $[u,t_1]$ is left-normed and $\deg [u,t_1]>\deg u$.
\epr

\begin{proposition}
\label{1-2}
The polynomials $[[t_1,t_2][t_3,t_4],t_5]$ and $[t_1,t_2][t_3,t_4][t_5,t_6]$ are identities for the algebra $F$.
\end{proposition}

\proof
Both polynomials are multilinear therefore it is enough to evaluate them on a spanning set of the algebra $F$. The algebra $F$ is spanned by elements of the type $C_1^nC_2^m u_1\cdots u_k$ where $n$, $m\ge 0$ and $u_i$ are left-normed commutators. Moreover if $k\ge 2$ then $u_1\cdots u_k$ is strongly central hence $C_1^nC_2^m u_1\cdots u_k$ will be central. But central elements vanish the commutators hence we consider substitutions by elements of the types $C_1^nC_2^m$ and $C_1^nC_2^m u$ only where $u$ is a left-normed commutator.

According to Lemmas~\ref{equalities}, \ref{leftright} one has $[C_1^nC_2^m, C_1^pC_2^q] = \sum_i w_i u_i$ where the $u_i$ are left-normed commutators, analogously for $[C_1^nC_2^mu, C_1^pC_2^q]$ and also for $[C_1^nC_2^m u, C_1^pC_2^q v]$. But then $[t_1,t_2][t_3,t_4]$ becomes $\sum_i w_iu_iv_i$ where $u_i$ and $v_i$ are left-normed commutators. The latter sum is (strongly) central and thus the first polynomial is an identity for $F$. The same procedure applied to the second polynomial yields a combination of products of three commutators  which is 0 in $F$ according to Lemma~\ref{3comm}.
\epr

\begin{corollary}
\label{2times3}
The identity $[t_1,t_2,t_5][t_3,t_4] + [t_1,t_2][t_3,t_4,t_5]=0$ holds in $F$.
\end{corollary}

\proof
It is another form of the first identity of Proposition~\ref{1-2}.
\epr

\noindent \textbf{Remark} The polynomial $[t_1,t_2][t_3,t_4]$ is a central polynomial for $U_2(M_{11}(E))$, and in particular $[t_1,t_2]^2$ is central as
well. On the other hand the latter polynomial is \textsl{not} central for $M_{11}(E)$. It is well known that for the matrix algebra $M_n(K)$ a polynomial
$f(t_1,\ldots,t_k)$ is central if and only if $f(A_1,\ldots,A_k)$ lies in the centre of the generic algebra $U_k(M_n(K))$ generated by $A_1$, \dots, $A_k$,
see \cite[Proposition 1.2, p. 171]{Procesi}. This is a sharp difference in the behaviour of the generic algebras for $M_n(K)$ and for $M_{ab}(E)$.

\medskip

\begin{proposition}
\label{s_4}
The standard polynomial $s_4= \sum (-1)^\sigma t_{\sigma(1)} t_{\sigma(2)} t_{\sigma(3)} t_{\sigma(4)}$ is an identity for $F$. Here $\sigma$ runs over the permutations of the symmetric group $S_4$, and $(-1)^\sigma$ stands for the sign of $\sigma$.
\end{proposition}

\proof
We write $s_4 = [t_1,t_2]\circ [t_3,t_4] - [t_1,t_3]\circ [t_2,t_4] + [t_1,t_4]\circ [t_2,t_3]$ where $a\circ b=ab+ba$.
As in Proposition~\ref{1-2} we shall substitute the variables by elements of the types $C_1^nC_2^m$ and $C_1^nC_2^m u$, $u$ a left-normed commutator.

First suppose $t_1=v_1u$, $t_i=v_i$ where $u$ is a left-normed commutator, and $v_i\in F$ are arbitrary. Then $[v_1u,v_2] = v_1[u,v_2] + [v_1,v_2]u$. The product of three commutators vanishes in $F$ hence
\[
[v_1u,v_2]\circ [v_3,v_4] = ([v_1,v_2]u) \circ [v_3,v_4]+ (v_1[u,v_2])\circ [v_3,v_4] = (v_1[u,v_2])\circ [v_3,v_4].
\]
Simple manipulations show that
\begin{eqnarray*}
(v_1[u,v_2])\circ [v_3,v_4] &=&
v_1[u,v_2][v_3,v_4] +[v_3,v_4]v_1[u,v_2]\\
&=& - v_1u[v_3,v_4,v_2] + v_1[v_3,v_4][u,v_2] + [v_3,v_4,v_1][u,v_2]\\
&=& - v_1u[v_3,v_4,v_2] - v_1[v_3,v_4,v_2]u - [v_3,v_4,v_1,v_2]u.
\end{eqnarray*}
Then one obtains by the Jacobi identity 
\begin{eqnarray*}
s_4(v_1u,v_2,v_3,v_4) &=&
-v_1u([v_3,v_4,v_2] - [v_2,v_4,v_4] + [v_2,v_3,v_4])\\
&&- v_1([v_3,v_4,v_2] - [v_2,v_4,v_4] + [v_2,v_3,v_4])u\\
&&-([v_3,v_4,v_1,v_2] - [v_2,v_4,v_1,v_3] + [v_2,v_3,v_1,v_4])u\\
&=& -([v_3,v_4,v_1,v_2] - [v_2,v_4,v_1,v_3] + [v_2,v_3,v_1,v_4])u.
\end{eqnarray*}
The latter sum equals, once again by Jacobi,
\[
([[v_1,v_2],[v_4,v_3]] + [[v_1,v_3],[v_2,v_4]] + [[v_1,v_4],[v_3,v_2]])u
\]
and this vanishes as a combination of products of three commutators.

Now we consider the substitution of $t_i$ by $C_1^{n_i}C_2^{m_i}$, $1\le i\le 4$. First one defines on $K[X;Y]$ an automorphism $'$ of order two by
$x_i\mapsto x_i'$, $y_i\mapsto y_i'$, and then extending to the whole supercommutative algebra. It is easy to see that for every $D=(d_{ij})\in F$ it holds
$d_{22}=d_{11}'$ and $d_{21} = d_{12}'$. This notation agrees also with the formula for the product $C_1^n C_2^m$ in (\ref{explicitproduct}). Suppose
$D_i=C_1^{n_i}C_2^{m_i}$. Hence we may consider $D_i  =
\begin{pmatrix} a_i&b_i\\ b_i' & a_i' \end{pmatrix}$, for some $a_i\in K[X;Y]_0$ and $b_i\in K[X;Y]_1$.

Put $D=(d_{ij})=s_4(D_1,D_2,D_3,D_4)\in F$; according to the above it suffices to prove that $d_{11}=d_{21}=0$.

Direct computation shows that the $(1,1)$-entry of $[D_1,D_2]\circ [D_3,D_4]$ is
\begin{eqnarray*}
&&2(b_1b_2'+b_1'b_2)(b_3b_4'+b_3'b_4)\\
&+& (b_1(a_2'-a_2)-b_2(a_1'-a_1))(b_4'(a_3'-a_3)-b_3'(a_4'-a_4))\\
&+&(b_3(a_4'-a_4)-b_4(a_3'-a_3))(b_2'(a_1'-a_1)-b_1'(a_2'-a_2)).
\end{eqnarray*}
Now writing down the analogous expressions for the remaining two  summands of $s_4$ as above, and summing up we have
\begin{eqnarray*}
d_{11}&=& 2(b_1b_2'+b_1'b_2)(b_3b_4'+b_3'b_4)-2(b_1b_3'+b_1'b_3)(b_2b_4'+b_2'b_4)\\
&+& 2(b_1b_4'+b_1'b_4)(b_2b_3'+b_2'b_3)\\
                                                 &+&(b_1(a_2'-a_2)-b_2(a_1'-a_1))(b_4'(a_3'-a_3)-b_3'(a_4'-a_4))\\
                                                 &+&(b_3(a_4'-a_4)-b_4(a_3'-a_3))(b_2'(a_1'-a_1)-b_1'(a_2'-a_2))\\
                                                 &-&(b_1(a_3'-a_3)-b_3(a_1'-a_1))(b_4'(a_2'-a_2)-b_2'(a_4'-a_4))\\
                                                 &-&(b_2(a_4'-a_4)-b_4(a_2'-a_2))(b_3'(a_1'-a_1)-b_1'(a_3'-a_3))\\
                                                 &+&(b_1(a_4'-a_4)-b_4(a_1'-a_1))(b_3'(a_2'-a_2)-b_2'(a_3'-a_3))\\
                                                 &+&(b_2(a_3'-a_3)-b_3(a_2'-a_2))(b_4'(a_1'-a_1)-b_1'(a_4'-a_4)).
\end{eqnarray*}
The last six expressions above cancel altogether, and we are left with
\[
d_{11}= 4(b_1'b_2'b_3b_4 - b_1'b_3'b_2b_4 +b_1'b_4'b_2b_3+ b_2'b_3'b_1b_4- b_2'b_4'b_1b_3+b_3'b_4'b_1b_2).
\]
As $D_i=C_1^{n_i}C_2^{m_i}$ it follows by (\ref{explicitproduct}) that $b_i=x_1^{n_i}Q_{m_i-1}y_2+x_2'^{m_i}q_{n_i-1}y_1+c_i$ for some $c_i\in
K[X;Y]^{(3)}$. Also $b_i'b_j'b_kb_l=g(i,j,k,l) y_1y_2y_1'y_2'$ where
\begin{eqnarray*}
g(i,j,k,l) &=&
x_1^{n_l}x_1'^{n_j} q_{n_k-1}q_{n_i-1} x_2^{m_i}x_2'^{m_k} Q_{m_j-1}Q_{m_l-1}\\
&+& x_1^{n_k}x_1'^{n_i} q_{n_j-1}q_{n_l-1} x_2^{m_j}x_2'^{m_l} Q_{m_i-1}Q_{m_k-1}\\
&-& x_1^{n_k}x_1'^{n_j} q_{n_l-1}q_{n_i-1} x_2^{m_i}x_2'^{m_l} Q_{m_j-1}Q_{m_k-1}\\
&-& x_1^{n_l}x_1'^{n_i} q_{n_j-1}q_{n_k-1} x_2^{m_j}x_2'^{m_k} Q_{m_i-1}Q_{m_l-1}.
\end{eqnarray*}
Thus for $d_{11}=s_4(D_1,D_2,D_3,D_4)_{11}$ we have
\begin{eqnarray*} d_{11}&=&
(g(1,2,3,4)-g(1,3,2,4)+g(1,4,2,3))y_1y_2y_1'y_2'\\
 &+& (g(2,3,1,4)-g(2,4,1,3)+g(3,4,1,2))y_1y_2y_1'y_2'.
\end{eqnarray*}
Expanding the sum of the $g(i,j,k,l)$ above we arrive at
\begin{eqnarray*}
& & Q_{m_2-1}Q_{m_3-1}q_{n_1-1}q_{n_4-1}(x_2^{m_4}x_2'^{m_1}-x^{m_1}x_2'^{m_4})(x_1^{n_3}x_1'^{n_2}-x_1^{n_2}x_1'^{n_3})\\
&+& Q_{m_1-1}Q_{m_4-1}q_{n_2-1}q_{n_3-1}(x_2^{m_3}x_2'^{m_2}-x_2^{m_2}x_2'^{m_3})(x_1^{n_4}x_1'^{n_1}-x_1^{n_1}x_1'^{n_4})\\
&+& Q_{m_2-1}Q_{m_4-1}q_{n_1-1}q_{n_3-1}(x_2^{m_1}x_2'^{m_3}-x_2^{m_3}x_2'^{m_1})(x_1^{n_4}x_1'^{n_2}-x_1^{n_2}x_1'^{n_4})\\
&+& Q_{m_1-1}Q_{m_3-1}q_{n_2-1}q_{n_4-1}(x_2^{m_2}x_2'^{m_4}-x_2^{m_4}x_2'^{m_2})(x_1^{n_3}x_1'^{n_1}-x_1^{n_1}x_1'^{n_3})\\
&+& Q_{m_1-1}Q_{m_2-1}q_{n_3-1}q_{n_4-1}(x_2^{m_4}x_2'^{m_3}-x_2^{m_3}x_2'^{m_4})(x_1^{n_2}x_1'^{n_1}-x_1^{n_1}x_1'^{n_2})\\
&+& Q_{m_3-1}Q_{m_4-1}q_{n_1-1}q_{n_2-1}(x_2^{m_1}x_2'^{m_2}-x_2^{m_2}x_2'^{m_1})(x_1^{n_3}x_1'^{n_4}-x_1^{n_4}x_1'^{n_3}).
\end{eqnarray*}
By the relations from Section 2: $(x_1^nx_1'^m - x_1^mx_1'^n) = (x_1-x_1')(q_nq_{m-1}-q_mq_{n-1})$ and
$(x_2^nx_2'^m-x_2^mx_2'^n)=(x_2-x_2')(Q_nQ_{m-1}-Q_mQ_{n-1})$ it follows $d_{11}=0$.

Now we prove that $s_4(D_1,D_2,D_3,D_4)_{21}=d_{21}=0$. The approach is similar to that of $d_{11}$. Computing the $(2,1)$-entry of $[D_1,D_2]\circ[D_3,D_4]$ we obtain
\begin{eqnarray*}
&& 2(b_2'(a_1'-a_1)-b_1'(a_2'-a_2))(b_3b_4'+b_3'b_4) \\
&+&  2(b_4'(a_3'-a_3)-b_3'(a_4'-a_4))(b_1b_2'+b_1'b_2)\\
&=& 2(a_1'-a_1) (b_2'b_3'b_4 + b_2'b_3b_4') - 2(a_2' - a_2) (b_1'b_3'b_4 + b_1'b_3b_4') \\
&+& 2(a_3'-a_3) (b_1b_2'b_4' + b_1'b_2b_4')- 2(a_4' - a_4)(b_1'b_2b_3' + b_1b_2'b_3').
\end{eqnarray*}
Permuting the indices we get that $(1/4)d_{21}$ equals
\begin{eqnarray*}
&& (a_1'-a_1) (b_2'b_3'b_4 - b_2'b_4'b_3 + b_3'b_4'b_2)
 - (a_2'-a_2) (b_1'b_3'b_4 - b_1'b_4'b_3 + b_3'b_4'b_1)\\
&& +(a_3'-a_3) (b_1'b_2'b_4 - b_1'b_4'b_2 + b_2'b_4'b_1) - (a_4'-a_4)(b_1'b_2'b_3 - b_1'b_3'b_2 + b_2'b_3'b_1).
\end{eqnarray*}
Substitute $b_i=x_1^{n_i}Q_{m_i-1}y_2+x_2'^{m_i}q_{n_i-1}y_1+c_i$, and $a_i=x_1^{n_i}x_2^{m_i}+d_i+e_i$ as in the previous case. Here $c_i\in
K[X;Y]^{(3)}$, $d_i\in K[X;Y]^{(2)}$, and $e_i\in K[X;Y]^{(4)}$.

One has $b_i'b_j'b_k = f_0(i,j,k)y_1y_1'y_2'+g_0(i,j,k)y_2y_1'y_2'$ with
\begin{eqnarray*}
f_0(i,j,k)&=& x_2'^{m_k} q_{n_k-1} (x_2^{m_i} Q_{m_j-1}x_1'^{n_j} q_{n_i-1} - x_2^{m_j}Q_{m_i-1} x_1'^{n_i}q_{n_j-1})\\
g_0(i,j,k)&=& x_1^{n_k} Q_{m_k-1}(x_2^{m_i}Q_{m_j-1}x_1'^{n_j} q_{n_i-1} - x_2^{m_j}Q_{m_i-1}x_1'^{n_i}q_{n_j-1}).
\end{eqnarray*}
By applying the equality $(x_2^n x_2'^m - x_2^m x_2'^n)=(x_2'-x_2) (Q_nQ_{m-1} - Q_mQ_{n-1})$ and after manipulations we get that
$f_0(i,j,k)-f_0(i,k,j)+f_0(j,k,i)$ equals
\begin{eqnarray*}
&& q_{n_k-1} q_{n_i-1} Q_{m_j-1}x_1'^{n_j}(x_2'^{m_k} x_2^{m_i} -x_2'^{m_i} x_2^{m_k})\\
&+& q_{n_k-1}q_{n_j-1}Q_{m_i-1}x_1'^{n_i}(x_2'^{m_j} x_2^{m_k} -x_2'^{m_k}x_2^{m_j})\\
&+& q_{n_j-1} q_{n_i-1} Q_{m_k-1}x_1'^{n_k}(x_2'^{m_i}x_2^{m_j} -x_2'^{m_j} x_2^{m_i})\\
&=& Q_{m_i}Q_{m_j-1} Q_{m_k-1} q_{n_i-1}(q_{n_k-1} x_1'^{n_j} -q_{n_j-1}x_1'^{n_k}) (x_2'-x_2)\\
&+& Q_{m_i-1} Q_{m_j} Q_{m_k-1} q_{n_j-1} (q_{n_i-1}x_1'^{n_k} -q_{n_k-1} x_1'^{n_i})(x_2' - x_2)\\
&+& Q_{m_i-1} Q_{m_j-1} Q_{m_k} q_{n_k-1} (q_{n_j-1}x_1'^{n_i} -q_{n_i-1}x_1'^{n_j})(x_2'-x_2).
\end{eqnarray*}
Thus $d_{21}=4f(x)y_1y_1'y_2'+ 4g(x)y_2y_1'y_2'$. We shall prove $f(x)=g(x)=0$. But
\begin{eqnarray*} (1/4)f(x) &=&
(x_1'^{n_1}x_2'^{m_1}-x_1^{n_1}x_2^{m_1})(f_0(2,3,4)-f_0(2,4,3)+f_0(3,4,2))\\
&-& (x_1'^{n_2}x_2'^{m_2} - x_1^{n_2}x_2^{m_2}) (f_0(1,3,4)-f_0(1,4,3) +f_0(3,4,1))\\
&+& (x_1'^{n_3}x_2'^{m_3}-x_1^{n_3}x_2^{m_3}) (f_0(1,2,4) -f_0(1,4,2) +f_0(2,4,1))\\
&-& (x_1'^{n_4}x_2'^{m_4}-x_1^{n_4}x_2^{m_4}) (f_0(1,2,3) -f_0(1,3,2) +f_0(2,3,1)).
\end{eqnarray*}
But $(x_1'^nx_2'^m - x_1^n x_2^m) = x_1'^n(x_2' -x_2) Q_{m-1} +x_2^m (x_1'-x_1)q_{n-1}$; substituting the $f_0$ by their defining equalities  in
$f_0(i,j,k)-f_0(i,k,j)+f_0(j,k,i)$ we get
$f(x)=(x_2'-x_2)^2f^{(1)}(x)+(x_2'-x_2)(x_1'-x_1)f^{(2)}(x)$.
Here
\begin{eqnarray*} f^{(1)}(x)&=& Q_{m_1-1} Q_{m_2} Q_{m_3-1} Q_{m_4-1} x_1'^{n_1}q_{n_2-1}(q_{n_4-1}x_1'^{n_3}-q_{n_3-1}x_1'^{n_4})\\
&+& Q_{m_1-1} Q_{m_2-1}Q_{m_3} Q_{m_4-1} x_1'^{n_1}q_{n_3-1} (q_{n_2-1} x_1'^{n_4}-q_{n_4-1}x_1'^{n_2})\\
&+& Q_{m_1-1} Q_{m_2-1} Q_{m_3-1} Q_{m_4} x_1'^{n_1} q_{n_4-1}(q_{n_3-1} x_1'^{n_2}-q_{n_2-1} x_1'^{n_3})\\
&-& Q_{m_1} Q_{m_2-1} Q_{m_3-1} Q_{m_4-1} x_1'^{n_2} q_{n_1-1}(q_{n_4-1} x_1'^{n_3}- q_{n_3-1} x_1'^{n_4})\\
&-& Q_{m_1-1} Q_{m_2-1} Q_{m_3} Q_{m_4-1} x_1'^{n_2} q_{n_3-1}(q_{n_1-1} x_1'^{n_4}- q_{n_4-1} x_1'^{n_1})\\
&-& Q_{m_1-1} Q_{m_2-1} Q_{m_3-1} Q_{m_4}x_1'^{n_2} q_{n_4-1}(q_{n_3-1} x_1'^{n_1}- q_{n_1-1}x_1'^{n_3})\\
&+& Q_{m_1} Q_{m_2-1} Q_{m_3-1} Q_{m_4-1} x_1'^{n_3} q_{n_1-1}(q_{n_4-1} x_1'^{n_2} - q_{n_2-1}x_1'^{n_4})\\
&+& Q_{m_1-1} Q_{m_2} Q_{m_3-1} Q_{m_4-1}x_1'^{n_3} q_{n_2-1}(q_{n_1-1} x_1'^{n_4}-q_{n_4-1} x_1'^{n_1})\\
&+& Q_{m_1-1} Q_{m_2-1} Q_{m_3-1} Q_{m_4}x_1'^{n_3} q_{n_4-1}(q_{n_2-1} x_1'^{n_1}-q_{n_1-1} x_1'^{n_2})\\
&-& Q_{m_1} Q_{m_2-1} Q_{m_3-1} Q_{m_4-1}x_1'^{n_4}  q_{n_1-1}(q_{n_3-1} x_1'^{n_2}- q_{n_2-1}x_1'^{n_3})\\
&-& Q_{m_1-1} Q_{m_2} Q_{m_3-1} Q_{m_4-1}x_1'^{n_4} q_{n_2-1}(q_{n_1-1} x_1'^{n_3}- q_{n_3-1}x_1'^{n_1})\\
&-& Q_{m_1-1}Q_{m_2-1} Q_{m_3-1}Q_{m_4} x_1'^{n_4} q_{n_3-1}(q_{n_2-1} x_1'^{n_1} -q_{n_1-1}x_1'^{n_2}).
\end{eqnarray*}
After simple manipulations we obtain $f^{(1)}(x)=0$. As for $f^{(2)}(x)$ we have
\begin{eqnarray*}
f^{(2)}(x)&=& q_{n_1-1}x_2^{m_1} Q_{m_2} Q_{m_3-1} Q_{m_4-1} q_{n_2-1}(q_{n_4-1}x_1'^{n_3}- q_{n_3-1}x_1'^{n_4}) \\
 &+& q_{n_1-1}x_2^{m_1} Q_{m_2-1} Q_{m_3} Q_{m_4-1} q_{n_3-1}(q_{n_2-1}x_1'^{n_4}- q_{n_4-1}x_1'^{n_2}) \\
&+&q_{n_1-1}x_2^{m_1}Q_{m_2-1} Q_{m_3-1} Q_{m_4}q_{n_4-1}(q_{n_3-1} x_1'^{n_2}- q_{n_2-1}x_1'^{n_3}) \\
&-&q_{n_2-1}x_2^{m_2} Q_{m_1}Q_{m_3-1} Q_{m_4-1}q_{n_1-1}(q_{n_4-1} x_1'^{n_3}-q_{n_3-1}x_1'^{n_4}) \\
&-&q_{n_2-1}x_2^{m_2}Q_{m_1-1} Q_{m_3}Q_{m_4-1} q_{n_3-1}(q_{n_1-1} x_1'^{n_4} -q_{n_4-1}x_1'^{n_1}) \\
&-&q_{n_2-1}x _2^{m_2}Q_{m_1-1}Q_{m_3-1} Q_{m_4}q_{n_4-1}(q_{n_3-1} x_1'^{n_1}-q_{n_1-1}x_1'^{n_3}) \\
&+&q_{n_3-1}x_2^{m_3}Q_{m_1} Q_{m_2-1}Q_{m_4-1} q_{n_1-1}(q_{n_4-1} x_1'^{n_2} - q_{n_2-1}x_1'^{n_4}) \\
&+&q_{n_3-1}x_2^{m_3}Q_{m_1-1} Q_{m_2}Q_{m_4-1}q_{n_2-1}(q_{n_1-1} x_1'^{n_4} -q_{n_4-1}x_1'^{n_1}) \\
&+&q_{n_3-1}x_2^{m_3}Q_{m_1-1} Q_{m_2-1}Q_{m_4} q_{n_4-1}(q_{n_2-1} x_1'^{n_1}- q_{n_1-1}x_1'^{n_2}) \\
&-&q_{n_4-1}x_2^{m_4}Q_{m_1}Q_{m_2-1} Q_{m_3-1} q_{n_1-1}(q_{n_3-1} x_1'^{n_2}- q_{n_2-1}x_1'^{n_3}) \\
&-&q_{n_4-1}x_2^{m_4}Q_{m_1-1}  Q_{m_2}Q_{m_4-1}q_{n_2-1}(q_{n_1-1} x_1'^{n_3}-q_{n_3-1}x_1'^{n_1}) \\
&-&q_{n_4-1}x_2^{m_4}Q_{m_1-1}Q_{m_2-1} Q_{m_4}q_{n_3-1}(q_{n_2-1} x_1'^{n_1} -q_{n_1-1}x_1'^{n_2})
\end{eqnarray*}
which in turn equals
\begin{eqnarray*}
&& q_{n_1-1}q_{n_2-1}q_{n_3-1} Q_{m_4-1}x_1'^{n_4} x_2^{m_1}(Q_{m_2-1} Q_{m_3}-Q_{m_2}Q_{m_3-1})\\
&+&q_{n_1-1}q_{n_2-1}q_{n_3-1}Q_{m_4-1} x_1'^{n_4}x_2^{m_2}(Q_{m_1} Q_{m_3-1}-Q_{m_1-1}Q_{m_3})\\
&+&q_{n_1-1}q_{n_2-1}q_{n_3-1}Q_{m_4-1} x_1'^{n_4} x_2^{m_3}(Q_{m_1-1} Q_{m_2}-Q_{m_1}Q_{m_2-1})\\
&-&q_{n_1-1}q_{n_2-1}q_{n_4-1}Q_{m_3-1} x_1'^{n_3}x_2^{m_1}(Q_{m_2-1} Q_{m_4}-Q_{m_2}Q_{m_4-1})\\
&-&q_{n_1-1}q_{n_2-1}q_{n_4-1}Q_{m_3-1} x_1'^{n_3}x_2^{m_2}(Q_{m_1} Q_{m_4-1}-Q_{m_1-1}Q_{m_4})\\
&-&q_{n_1-1}q_{n_2-1}q_{n_4-1}Q_{m_3-1} x_1'^{n_3}x_2^{m_4}(Q_{m_1-1} Q_{m_2}-Q_{m_1}Q_{m_2-1})\\
&+&q_{n_1-1}q_{n_3-1}q_{n_4-1}Q_{m_2-1} x_1'^{n_2}x_2^{m_1}(Q_{m_3-1} Q_{m_4}-Q_{m_3}Q_{m_4-1})\\
&+&q_{n_1-1}q_{n_3-1}q_{n_4-1} Q_{m_2-1}x_1'^{n_2}x_2^{m_3}(Q_{m_1} Q_{m_4-1}-Q_{m_1-1}Q_{m_4})\\
&+&q_{n_1-1}q_{n_3-1}q_{n_4-1} Q_{m_2-1}x_1'^{n_2}x_2^{m_4}(Q_{m_1-1} Q_{m_3}-Q_{m_1}Q_{m_3-1})\\
&-&q_{n_2-1}q_{n_3-1}q_{n_4-1}Q_{m_1-1} x_1'^{n_1}x_2^{m_2}(Q_{m_3-1} Q_{m_4}-Q_{m_3}Q_{m_4-1})\\
&-&q_{n_2-1}q_{n_3-1}q_{n_4-1} Q_{m_1-1}x_1'^{n_1}x_2^{m_3}(Q_{m_2} Q_{m_4-1}-Q_{m_2-1}Q_{m_4})\\
&-&q_{n_2-1}q_{n_3-1}q_{n_4-1} Q_{m_1-1}x_1'^{n_1}x_2^{m_4}(Q_{m_2-1}Q_{m_3}-Q_{m_2}Q_{m_3-1}).
\end{eqnarray*}
Applying $(x_2^n x_2'^m- x_2^mx_2'^n) = (x_2- x_2')(Q_nQ_ {m-1} -Q_mQ_{n-1})$ one sees that the first three summands in the final expression cancel out.
Repeat the procedure for the remaining three groups of 3 summands in each thus getting $f^{(2)}(x)=0$.

In order to show that $g(x)=0$, we observe that
\begin{eqnarray*} (1/4)g(x) &=&
    (x_1'^{n_1}x_2'^{m_1}-x_1^{n_1}x_2^{m_1}) (g_0(2,3,4)-g_0(2,4,3)+g_0(3,4,2))\\
&-& (x_1'^{n_2}x_2'^{m_2}-x_1^{n_2}x_2^{m_2}) (g_0(1,3,4)-g_0(1,4,3)+g_0(3,4,1))\\
&+& (x_1'^{n_3}x_2'^{m_3}-x_1^{n_3}x_2^{m_3}) (g_0(1,2,4)-g_0(1,4,2)+g_0(2,4,1))\\
&-& (x_1'^{n_4}x_2'^{m_4}-x_1^{n_4}x_2^{m_4}) (g_0(1,2,3)-g_0(1,3,2)+g_0(2,3,1)).
\end{eqnarray*}

Now consider the automorphism $\phi$ of $K[x_1,x_2,x_1',x_2']$ defined by $x_1\mapsto x_2'$, $x_2\mapsto x_1'$, $x_1'\mapsto x_2$ and $x_2'\mapsto x_1$.
Let us compute $\phi(f(x))$.

First observe that for each $i$, $\phi(x_1'^{n_i}x_2'^{m_i}-x_1^{n_i}x_2^{m_i})=-(x_1'^{m_i}x_2'^{n_i}-x_1^{m_i}x_2^{n_i})$ and for each $i$, $j$ and $k$,
\[
\phi(f_0(i,j,k))=-x_1^{m_k}Q_{n_k-1}(x_1'^{m_j}q_{m_i-1}x_2^{n_i}Q_{n_j-1}-x_1'^{m_i}q_{m_j-1}x_2^{n_j}Q_{n_i-1}).
\]
By the above observations, we can see that $-\phi(f(x))$ is equal to the expression obtained by $g(x)$ by permuting $m_i$ and $n_i$ for each $i$. Since
$f(x)=0$, for arbitrary integers $n_i$, $m_i$, $i\in \{1,2,3,4\}$, it follows that $g(x)=0$. Thus $s_4$ is an identity for $F$. \epr

\noindent \textbf{Remark 1} If one tries to prove the above Proposition by splitting $F$ replacing it with the algebra generated by
$\{x_ie_{11},x_i'e_{22},y_ie_{12},y_i'e_{21},i=1,2\}$, the resulting larger algebra does not satisfy $s_4$, since
\[s_4(y_1e_{12},y_1'e_{21},y_2e_{12},y_2'e_{21})=4y_1y_1'y_2y_2'(e_{11}+e_{22})\neq 0.\]

\noindent \textbf{Remark 2} Neither of the identities from Propositions~(\ref{1-2}), (\ref{s_4}) is a polynomial identity for $M_{11}(E)$. Clearly
$T(M_{11}(E))\subseteq T(F)$; it can be easily shown that all identities for $M_{11}(E)$ follow from $[[t_1,t_2][t_3,t_4],t_5]$ (though the latter is not
an identity for $M_{11}(E)$).

\section{The identities of $M_{11}(E)$}

Here we shall prove that the polynomials from Propositions~(\ref{1-2}), (\ref{s_4}) generate the T-ideal of $F$. To this end we make use of the results of  Popov \cite{popov}. Namely we shall need not only the concrete form of the basis of $T(M_{11}(E))$ but also the structure of the corresponding relatively free algebra. In order to make our exposition more self-contained we recall here the results from \cite{popov}.

Let $\Gamma_n$ be the vector space of the proper multilinear polynomials in $t_1$, \dots, $t_n$ in the free associative algebra $K\langle T\rangle$. We work with unitary algebras over a field of characteristic 0 hence the T-ideal of any algebra $A$ is generated by the intersections $T(A)\cap \Gamma_n$. The vector space $\Gamma_n$ is a left module over the symmetric group $S_n$. The action of $S_n$ on $\Gamma_n$ is by permuting the variables. We refer the reader to \cite[Chapter 12]{drenskybook} for all necessary information concerning the representations of $S_n$ and their applications to PI theory.

Let $W$ be the variety of algebras defined by $M_{11}(E)$ and put $\Gamma_n(W) = \Gamma_n/(\Gamma_n\cap T(W))$. Then clearly $\Gamma_n(W)$ is an
$S_n$-module (with the induced action). The result from \cite{popov} we need is the decomposition of $\Gamma_n(W)$ into irreducible submodules. It is well
known that the irreducibles for $S_n$ are described in terms of partitions of $n$ and Young tableaux. We refer once again to \cite{drenskybook} for that
description. Recall that the polynomial representations of the general linear group $GL_m$ are also described in terms of partitions (of not more than $m$
parts) and Young diagrams, see \cite{drenskybook}. Sometimes it is convenient to use the ''symmetrized'' version of a generator of an irreducible
$S_n$-module; it generates an irreducible $GL_m$-module, and vice versa via linearisation. Recall that one can linearise and go back to the symmetrized version of a  polynomial as char$K=0$.

The following polynomials were defined in \cite{popov}. Let $p\ge q\ge 2$ and $s\ge 0$. Set $\varphi_p^{(s)}=\varphi_p^{(s)}(t_1, \dots, t_p)$ and $\varphi_{p,q}^{(s)} = \varphi_{p,q}^{(s)}(t_1, \dots,t_p)$ as follows.

\begin{eqnarray*}
\varphi_p^{(s)}&=&\begin{cases}  \sum_{\sigma\in S_p}(-1)^{\sigma}[t_{\sigma(1)}, t_{\sigma(2)}] \dots [t_{\sigma(p)}, x_1^{(s)}] & p \text{ odd,}  \\
\sum_{\sigma\in S_p}(-1)^{\sigma}[t_{\sigma(1)},t_{\sigma(2)}]\dots[t_{\sigma(p-1)},t_{\sigma(p)}, t_1^{(s)}] & p \text{ even,} \\
\end{cases}\\
\varphi_{p,q}^{(s)}&=&
\begin{cases}
\sum_{\tau\in S_q}(-1)^{\tau}[t_{\tau(1)},t_{\tau(2)}]\dots[t_{\tau(q-1)},t_{\tau(q)}] \varphi_p^{(s)} & q\text{ even,} \\
\sum_{\sigma\in S_p}(-1)^{\sigma}[t_{\sigma(1)},t_{\sigma(2)}] \dots[t_{\sigma(p-1)},t_{\sigma(p)}] \varphi_q^{(s)} & p \text{ even, } q \text{ odd} \end{cases}
\end{eqnarray*}
and $\varphi_{p,q}^{(s)}=\sum (-1)^{\sigma\tau}[x_{\tau(1)},x_{\tau(2)}] \cdots[x_{\sigma(1)},x_{\sigma(2)}] \cdots
[x_{\sigma(p)}, x_1^{(s)}, x_{\tau(q)}]$  when both $p$ and $q$ are odd.  Here the summation runs over $\sigma\in S_p$, $\tau\in S_q$, and $[a,b^{(s)}]$ stands for $[a,b,\ldots,b]$ with $s$ entries of $b$.

We observe that for $s=0$, some of these polynomials are not defined. For example, for $p$ even, $\varphi_p^{(0)}(x_1,\dots,x_p)$ is the standard
polynomial of degree $p$, while for $p$ odd, it is not defined, since the associated diagram is a single column, which is associated to the standard
polynomial of degree $p$ that is not a proper polynomial when $p$ is odd. In a similar way $\varphi_{p,q}^{(0)}=s_q(x_1,\dots,x_q)s_p(x_1,\dots,x_p)$, if
$p$ and $q$ are even and $\varphi_{p,q}^{(0)}=\sum (-1)^{\sigma\tau}[x_{\tau(1)},x_{\tau(2)}] \cdots[x_{\sigma(1)},x_{\sigma(2)}] \cdots
[x_{\sigma(p)},x_{\tau(q)}]$ if $p$ and $q$ are odd, and it is not defined for the remaining cases.

Set  $M_p^{(s)}$ to be the $S_n$-submodule of  $\Gamma_n(W)$ generated by $\varphi_p^{(s)}$,  $n=p+s$, and $M_{p,q}^{(s)}$ the
$S_n$-submodule of $\Gamma_n(W)$ generated by $\varphi_{p,q}$,  $n=p+q+s$.

Let $f(t_1,\ldots,t_n)\in \Gamma_n$ be a proper multilinear polynomial and suppose $d=D_\alpha$ is a Young tableau associated to a diagram $D$ of a partition of $n$. That is we fill the boxes of the diagram $D$, along the rows, with the numbers of the permutation $\alpha$ of $n$. Form the Young semi-idempotent $e(d)$ of $d$ and denote by $M(d,f)$ the $S_n$-module generated by $e(d) f$. It is well known that it is either 0 or  irreducible.

The polynomials $\varphi_p^{(s)}$ and $\varphi_{p,q}^{(s)}$ are not multilinear but we consider their complete linearisations. The description of
$\Gamma_n(W)$ given in \cite{popov} is based on the following results.
\begin{enumerate}
\item
If $p\ge 2$, $s\ge 0$ then $\varphi_p^{(s)}$ and $\varphi_{p,q}^{(s)}$ are not polynomial identities for $M_{11}(E)$.
\item
Let $n=p+s$ and let $D=(s+1,1^{p-1})$ be a partition of $n$ with an associated  Young tableau $d=D_\alpha$. If $f(t_1,\ldots,t_n)\in \Gamma_n(W)$ then $M(d,f)\subseteq M_p^{(s)}$
\item
Let $n=p+q+s$, $d=D_\alpha$ with $D=(s+2,2^{q-1},1^{p-q})$ and $f=f(x_1,\dots,x_n)\in \Gamma_n(W)$. Then $M(d,f)\subseteq  M_{p,q}^{(s)}$.
\item
If the second row of a diagram $D$ contains at least 3 boxes then for every $f$ the module $M(d,f)$ is 0 in $\Gamma_n(W)$.
\item
The decomposition of $\Gamma_n(W)$ in irreducibles is
\[
\Gamma_n(W)=(\oplus_{p+s=n}M_{p}^{(s)}) \bigoplus (\oplus_{p+q+s=n} M_{p,q}^{(s)}).
\]
\end{enumerate}

Popov deduced, as a corollary to the above listed results, the main theorem of \cite{popov}, namely that the T-ideal of $M_{11}(E)$ is generated by the polynomials $[[t_1,t_2]^2, t_1]$ and $[t_1,t_2,[t_3,t_4],t_5]$.

\section{The identities of $F$}

\begin{theorem}
\label{mainthm}
Let $K$ be a field of characteristic 0. The polynomials
\begin{equation}
\label{basisF}
[[t_1,t_2][t_3,t_4],t_5], \quad [t_1,t_2] [t_3,t_4] [t_5,t_6], \quad s_4(t_1,t_2,t_3,t_4)
\end{equation}
form a basis of the polynomial identities for the algebra $F=K[C_1,C_2]$.
\end{theorem}

\proof
First we shall prove that the polynomials $\varphi_4^{(s)}$ lie in the T-ideal generated by the three polynomials from our theorem. Recall that
\[
\varphi_p^{(s)}=\begin{cases}  \sum_{\sigma\in S_p}(-1)^{\sigma}[t_{\sigma(1)}, t_{\sigma(2)}] \dots [t_{\sigma(p)}, x_1^{(s)}] & p \text{ odd,}  \\
\sum_{\sigma\in S_p}(-1)^{\sigma}[t_{\sigma(1)},t_{\sigma(2)}]\dots[t_{\sigma(p-1)},t_{\sigma(p)}, t_1^{(s)}] & p \text{ even.}
\end{cases}
\]
Also it is immediate that $\varphi_4^{(0)}(t_1,t_2,t_3,t_4) = 4s_4(t_1,t_2,t_3,t_4)$. Suppose $s\ge 1$. By  the identity from Corollary~\ref{2times3} we
obtain that $(1/4) \varphi_p^{(s)}(t_1,t_2,t_3,t_4)$ equals
\begin{eqnarray*}
&&[t_1,t_2][t_3,t_4,t_1^{(s)}] + (-1)^{s}[t_3,t_4,t_1^{(s)}] [t_1,t_2] + [t_2,t_3] [t_1,t_4,t_1^{(s)}] \\
&+& (-1)^{s}[t_1,t_4,t_1^{(s)}][t_2,t_3] - [t_2,t_4]  [t_1, t_3, t_1^{(s)}] - (-1)^{s}[t_1, t_3, t_1^{(s)}][t_2,t_4].
\end{eqnarray*}
Thus $(1/4) \varphi_p^{(s)}(t_1,t_1t_2,t_3,t_4)$ will be equal to
\begin{eqnarray*}
&& [t_1,t_3] t_2 [t_1,t_4,t_1^{(s)}] + [t_1,t_4] [t_1,t_3, t_1^{(s)}] t_2 - [t_1,t_4] t_2 [t_1,t_3,t_1^{(s)}] \\
&-& [t_1,t_3][t_1,t_4, t_1^{(s)}] t_2 - [t_1,t_4] [t_2, t_3,  t_1^{(s+1)}] - [x_3,x_4] [x_1,x_2, x_1^{(s+1)}] \\                                                                              &+& [t_1,t_3] [t_2,t_4,t_1^{(s+1)}]
+ (1/4)t_1\varphi_4^{(s)}(t_1, t_2,t_3, t_4).
\end{eqnarray*}
Analogously for $(1/4) \varphi_p^{(s)}(t_1,t_2t_1,t_3,t_4)$ we obtain
\begin{eqnarray*}
&&t_2[t_1,t_3][t_1,t_4,t_1^{(s)}] + (-1)^s [t_1,t_4, t_1^{(s)}] t_2[t_1,t_3] - t_2[t_1,t_4][t_1, t_3, t_1^{(s)}]\\
&-&(-1)^s[t_1, t_3, t_1^{(s)}] t_2[t_1,t_4] - [t_1,t_2][t_3,t_4, t_1^{(s+1)}] - [t_2,t_3][t_1,t_4,t_1^{(s+1)}]\\
&+& [t_2,t_4] [t_1,t_3,t_1^{(s+1)}] + (1/4) \varphi_4^{(s)} (t_1, t_2,t_3, t_4) t_1.
\end{eqnarray*}
Therefore for $(1/4) ( \varphi_p^{(s)}(t_1,t_1t_2,t_3,t_4) + \varphi_p^{(s)}(t_1,t_2t_1,t_3,t_4))$ we have
\begin{eqnarray*}
& & (1/4)(t_1\circ \varphi_4^{(s)}(t_1,t_2, t_3,t_4)- \varphi_4^{(s+1)} +  [t_1,t_3] t_2[t_1,t_4, t_1^{(s)}]\\
&-& (-1)^s[t_1,t_3,t_1^{(s)}]t_2[t_1,t_4] - [t_1,t_4]t_2 [t_1,t_3, t_1^{(s)}] +(-1)^s[t_1,t_4,t_1^{(s)}]t_2[t_1,t_3].
\end{eqnarray*}
An easy computation shows that the first two identities from the theorem, together with the fact that commutation by a fixed element is a derivation, imply the identity $[[t_1,t_2]t_3[t_4,t_5], t_6]=0$. Using the fact that the product of 3 commutators is 0 we obtain
$[t_1,t_2]t_3[t_4,t_5,t_6]+ [t_1,t_2,t_6]t_3[t_4,t_5]=0$. Repeating several times we arrive at
\begin{eqnarray*}
[t_1,t_3]t_2[t_1,t_4,t_1^{(s)}]- (-1)^s[t_1,t_3,t_1^{(s)}] t_2[t_1,t_4] &=&0,\\
{}[t_1, t_4] t_2[t_1,t_3,t_1^{(s)}]- (-1)^s[t_1,t_4,t_1^{(s)}] t_2[t_1,t_3]&=&0.
\end{eqnarray*}
In this way $\varphi_4^{(s+1)}(t_1,t_2,t_3,t_4)$ equals
\[\varphi_4^{(s)}(t_1,t_2,t_3,t_4)\circ t_1-\varphi_4^{(s)}(t_1,t_1t_2,t_3,t_4)-\varphi_4^{(s)}(t_1,t_2t_1,t_3,t_4),\]
and we may proceed by induction on $s$ since $\varphi_4^{(0)} = 4s_4$.

\medskip

As we work with unitary algebras and char$K=0$ we can consider only the multilinear proper identities. Denote by $I$ the T-ideal generated by the identities from the theorem and let $V$ be the variety of unitary algebras determined by $I$; we shall study the $S_n$-module $\Gamma_n(V) = \Gamma_n/\Gamma_n\cap I$. But the T-ideal of $M_{11}(E)$ is contained in $I$ hence $\Gamma_n(V)$ is a homomorphic image of $\Gamma_n(W)$, and we have to determine which of the irreducibles in $\Gamma_n(W)$ vanish modulo the T-ideal $I$.

But it is easy to see that the polynomials $\varphi_p^{(s)}$, $p\ge 5$, are products of at least three commutators, and as such they follow from the second identity of the theorem. The same holds for $\varphi_{p,q}^{(s)}$ whenever $p+q\ge 5$. Also we showed above that $\varphi_4^{(s)}$ lies in $I$ for every $s\ge 0$. Therefore
\[
\Gamma_n(V) = M_2^{(n-2)} \oplus M_3^{(n-3)} \oplus M_{2,2}^{(n-4)}.
\]
Hence in order to complete the proof it suffices to see that the generators of the irreducible modules from the above decomposition are all non-zero modulo $I$. But $\varphi_2^{(s)} (t_1,t_2) = 2[t_1,t_2, t_1^{(s)}]$ and $\varphi_2^{(s)}(C_1,C_2)\ne 0$ according to \cite[Lemma 9]{pktcm}.

Analogously $\varphi_3^{(s)} (t_1,t_2,t_3) = 2([t_1,t_2][t_3,t_1^{(s)}] - [t_1,t_3] [t_2,t_1^{(s)}])$ and we obtain
$\varphi_3^{(s)}(C_1,C_2,[C_1,C_2])=-4[C_1,C_2][C_2,C_1^{(s+1)}]\ne 0$ due to \cite[Lemma 9]{pktcm}.

Finally $\varphi_{2,2}^{(s)}=4[t_1,t_2][t_1,t_2,t_1^{(s)}]$ and $\varphi_{2,2}^{(s)}(C_1,C_2)=4[C_1,C_2][C_1,C_2,C_1^{(s)}]$ is non-zero for the same
reason as above. \epr

In \cite{gordienko}, A. Gordienko studied the identities satisfied by the algebra $A_1\subseteq UT_3(K)$ where $UT_3(K)$ are the upper triangular matrices
of order 3 over $K$. The algebra $A_1$ consists of all matrices whose $(1,1)$ and $(3,3)$ entries are equal. He deduced that $T(A_1)$ is generated by the
same three identities as in our Theorem~\ref{mainthm}. Moreover Gordienko described basis of the vector space $P_n(A_1)$ and $\Gamma_n(A_1)$ of the
multilinear and the proper multilinear elements of degree $n$ modulo the identities of $A_1$, respectively. Combining our result with Gordienko's theorem
we obtain the following corollary.

\begin{corollary}
The algebras $F = K[C_1,C_2]$ and $A_1$ are PI equivalent over a field of characteristic 0.
\end{corollary}

It will be interesting to know whether these two algebras remain PI equivalent if the field $K$ is infinite and of characteristic $p>2$.

\medskip

\section{A description of the subvarieties}

In this section we shall describe the subvarieties of the variety $V$ of unitary algebras defined by the identities of the algebra $F$. We shall work over a field of characteristic 0. Recall that $\Gamma_n(V) =  M_2^{(n-2)} \oplus M_3^{(n-3)} \oplus M_{2,2}^{(n-4)}$. Therefore we have to find the consequences of degree $n+1$ in $\Gamma_{n+1}(V)$ of the polynomials $\varphi_2^{(n-2)}$, $\varphi_{2,2}^{(n-4)}$, $\varphi_3^{(n-3)}$.

\begin{proposition}
\label{conseq}
The consequences of degree $n+1$  are as follows.

(1) The polynomials $\varphi_2^{(n-1)}$, $\varphi_{2,2}^{(n-3)}$, and  $\varphi_3^{(n-2)}$, follow from $\varphi_2^{(n-2)}$.

(2) The polynomials $\varphi_{2,2}^{(n-3)}$ and $\varphi_3^{(n-2)}$ follow from $\varphi_{2,2}^{(n-4)}$.

(3) The polynomials $\varphi_{2,2}^{(n-3)}$ and $\varphi_3^{(n-2)}$  follow from $\varphi_3^{(n-3)}$.
\end{proposition}

\proof Let $u_2^{(n-2)}(t_1,t_2,t_3)$ be the multihomogeneous component of the polynomial $\varphi_2^{(n-2)}(t_1+t_3,t_2)$ that is linear in $t_3$, and let
$u_{2,2}^{(n-4)}(t_1,t_2,t_3)$ be the component that is linear in $t_2$ and in $t_3$ of $\varphi_{2,2}^{(n-4)}(t_1,t_2+t_3)$. Clearly the identities
$u_2^{(n-2)}$ and $\varphi_2^{(n-2)}$ are equivalent, and also $u_{2,2}^{(n-4)}$ and $\varphi_{2,2}^{(n-4)}$ are equivalent, as the former are partial
linearisations of the latter. Then
\begin{eqnarray*}
u_2^{(n-2)} &=& -2([t_2,t_3,t_1^{(n-2)}]+(n-3)[t_2,t_1,t_3,t_1^{(n-3)}]+[t_2,t_1^{(n-2)},t_3]);\\
u_{2,2}^{(n-4)} &=& -4([t_1,t_2][t_3,t_1^{(n-3)}]+[t_1, t_3] [t_2, t_1^{(n-3)}]).
\end{eqnarray*}

In order to prove the proposition we verify which of the generators of the irreducible modules of $\Gamma_{n+1}(M_{1,1})$ are consequences of each
generator of the irreducible modules of $\Gamma_n(M_{1,1})$.

1. We compute directly that $\varphi_2^{(n-1)} = \varphi_2^{(n-2)}t_1 - t_1\varphi_2^{(n-2)}$. Similarly
\begin{eqnarray*}
\varphi_3^{(n-2)}(t_1,t_2,t_3) &= & (-1)^n (u_2^{(n-2)}(t_1, t_3, t_1t_2) - u_2^{(n-2)}(t_1, t_2,t_1t_3) \\
&+& t_1(u_2^{(n-2)} (t_1,t_2,t_3)- u_2^{(n-2)}(t_1, t_3,  t_2))/(n-2)\\
&+&  (n-1)(\varphi_2^{(n-2)}(t_1,t_2)t_3- \varphi_2^{(n-2)}(t_1, t_3) t_2)).
\end{eqnarray*}
Also $\varphi_{2,2}^{(n-3)} =  u_2^{(n-2)}(t_1,t_1t_2, t_2)  +2\varphi_2^{(n-2)}t_2-t_1u_2^{(n-2)}(t_1,t_2,t_2)$ when $n$ is odd, and similarly $\varphi_{2,2}^{(n-3)} = (u_2^{(n-1)}(t_1, t_2, t_2)  - [\varphi_2^{(n-2)}(t_1,t_2), t_2])/(n-2) - u_2^{(n-2)}(t_1, t_2,[t_1,t_2])$ if $n$ is even.

2. One computes directly that
\begin{eqnarray*}
\varphi_{2,2}^{(n-3)} & = & t_1 u_{2,2}^{(n-4)} (t_1,t_2,t_2)- u_{2,2}^{(n-4)}(t_1,t_1t_2, t_2);\\
\varphi_3^{(n-3)} & = & (u_{2,2}^{(n-4)}(t_1,t_2,t_1t_3) - u_{2, 2}^{(n-4)} (t_1,t_1t_2,t_3))/2.
\end{eqnarray*}
Moreover $\varphi_2^{(n-1)}$ is not a consequence of $\varphi_{2,2}^{(n-4)}$. In order to see this observe that $\varphi_{2,2}^{(n-4)}(C_1,C_2)$ is strongly central in $F$. Hence all its consequences will be central but $\varphi_2^{(n-2)}(C_1,C_2)$ is not central.

3. Finally, for $\varphi_{2,2}^{n-4}$, one checks directly that
\begin{eqnarray*}
\varphi_3^{(n-2)}& = & t_1 \varphi_3^{(n-3)}+ \varphi_3^{(n-3)} t_1 -\varphi_3^{(n-3)}(t_1,t_1t_2,t_3)- \varphi_3^{(n-3)}(t_1, t_2t_1, t_3);\\
\varphi_{2,2}^{(n-3)} & = & \varphi_3^{(n-3)}(t_1,t_2,[t_1,t_2]).
\end{eqnarray*}
As in (2), $\varphi_2^{(n-1)}$ is not a consequence of  $\varphi_3^{(n-3)}$.
\epr

\medskip

Let $I$ and $J$ be two T-ideals in $K\langle X\rangle$. Following \cite{Kemer2}, we say that $I$ and $J$ are \textsl{asymptotically equivalent} if for all
sufficiently large $n$ it holds $I\cap B_n=J\cap B_n$, where $B_n$ states for the proper polynomials of degree $n$. As we consider unitary algebras over a field of characteristic 0 we can
substitute $B_n$ by $\Gamma_n$.

\begin{corollary}
\label{asympt} If $U$ is a proper subvariety of $V$ then $U$ is asymptotically equivalent to either var$(K)$ or to var$(UT_2(K))$.
\end{corollary}

\proof
The variety $U$ satisfies, for some $n$, at least one of the identities $\varphi_2^{(n-2)}(x_1,x_2)$, $\varphi_3^{(n-3)}(x_1,x_2,x_3)$,
$\varphi_{2,2}^{(n-4)}(x_2,x-2)$.

If it satisfies some $\varphi_2^{(n-2)}$ then by Proposition~\ref{conseq} every commutator of sufficiently large degree will vanish on $U$. Therefore $U$ is asymptotically equivalent to var$(K)$. (Recall the latter is generated by the identity $[t_1,t_2]$.)

If on the other hand $U$ satisfies $\varphi_{2,2}^{(n-4)}$ or  $\varphi_3^{(n-3)}$ then every proper polynomial of sufficiently large degree, that is a product of two commutators, will vanish on $U$. But the T-ideal of $UT_2(K)$ is generated by the product of two commutators, hence $V$ is asymptotically equivalent to var$(UT_2(K))$.
\epr

\section{Identities in two variables for $M_{11}(E)$}

In this section we apply some of the results obtained above in order to describe the identities in two variables for the algebra $M_{11}(E)$ over a field of characteristic 0. Recall that in this way we determine the identities in two variables for $E\otimes E$ as well. The identities in two variables for $M_2(K)$, char$K=0$ were described by Nikolaev in \cite{nikolaev}. We shall employ some of the results obtained  in \cite{nikolaev}.

A word about the notation we shall use here. In order not to accumulate indices we shall use the letters $x$ and $y$ for the variables.

The main theorem in \cite{nikolaev} states that the identities in two variables for $M_2(K)$ all follow from the Hall polynomial $h(x,y) = [[x,y]^2,x]$. Recall that $T(M_2(K))$ is generated by $h$ and by $s_4$; one usually writes $h(x,y,z) = [[x,y]^2,z]$. This form of $h$ is equivalent to the one above \textsl{modulo} the standard polynomial $s_4$, see for example \cite{jcpk}; otherwise the two forms of the Hall polynomial are not equivalent.

Let $H$ be the variety of (unitary) algebras determined by the polynomial $h$. As in the previous section we shall work with the proper multilinear elements only. Define
\[
f_{kmn} = [x,y]^k [x,y,x^{(m)}, y^{(n)}]; \qquad d_{k,l}(x,y) = [x,y]^k [x,y,x^{(l)}].
\]
Nikolaev proved that the polynomials $f_{kmn}$, $k$, $m$, $n\ge 0$ span the vector space of all proper multilinear polynomials modulo the identity $h$.
Moreover he proved that $\Gamma_n(H)=\oplus M_{k,l}$ where the sum is over all $k$, $l\ge 0$ such that $2k+l+2=n$. The irreducible $S_n$-modules $M_{k,l}$
are generated by the complete linearisations of the polynomials $d_{k,l}$ given above.

But $h(x,y)=[[x,y]^2,x]$ vanishes on $M_{1,1}(E)$. Therefore $\Gamma_n(M_{11}(E))$ is a homomorphic image of $\Gamma_n(H)$. We denote by $J$ the T-ideal of the identities in two variables for $M_{11}(E)$, thus we will have a decomposition $\Gamma_n(M_{1,1})=\oplus M_{k,l}\pmod{J}$. Here the sum is taken over (some of the) $k$ and $l$. Thus we have to check which of the polynomials $d_{k,l}$ are identities for $M_{11}(E)$ and which are not.

\begin{lemma}
The polynomial $d_{k,l}$ is not an identity for $M_{1,1}(E)$ when $k < 2$.
\end{lemma}

\proof Suppose $k=0$, then $d_{0,l}=[x,y,x^{(l)}]$, and $d_{0,l}(C_1,C_2)\ne 0$. Therefore $d_{0,l}$ is not an identity for $M_{11}(E)$. If $k=1$ we have
$d_{1,l}=[x,y][x,y,x^{(l)}]$ and as above $d_{1,l}$ is not an identity for $M_{11}(E)$. \epr

\begin{lemma}
\label{dkl}
The polynomials $d_{k,l}$, $k\ge 2$, are identities for $M_{1,1}(E)$. If  $k\ge 3$ then $d_{k,l}$ follows from $d_{2,l}$.
\end{lemma}

\proof
Suppose $k\ge 2$, then for every $l\ge 0$ the polynomial $d_{k,l}$ is a product of three commutators. Therefore $d_{k,l}(C_1,C_2)=0$ and consequently $d_{k,l}\in T(M_{11}(E))$. The second statement is immediate.
\epr

Thus we have the following corollary.

\begin{corollary}
\label{manypi} All identities in two variables for $M_{11}(E)$ are consequences of $[[x,y]^2,x]$ and the polynomials $d_{2,l} = [x,y]^2[x,y,x^{(l)}]$,
$l\ge 0$.
\end{corollary}

\begin{theorem}
\label{twovar}
All identities in two variables for $M_{11}(E)$ follow from the two polynomials $h=[[x,y]^2,x]$ and $d=[x,y]^3$.
\end{theorem}

\proof
Clearly $h$ and $d$ are identities for $M_{11}(E)$. In order to prove the theorem it suffices to show, according to Lemma \ref{dkl}, that all polynomials $d_{2,l} = [x,y]^2 [x,y,x^{(l)}]$, are consequences of $d$ and $h$. We shall induct on $l$. The base of the induction is $l=0$ when $d_{2,0}=d$. Write then
\begin{eqnarray*}
[x,y]^2[x,y,x^{(l)}] &=& [x,y]^2([x,y,x^{(l-1)}]x -x[x,y,x^{(l-1)}])\\ &=&[x,y]^2[x,y,x^{(l-1)}]x-[x,y]^2x[x,y,x^{(l-1)}].
\end{eqnarray*}
As $[x,y]^2$ commutes with $x$ (and with $y$) we have that $d_{2,l}$ follows from $d_{2,l-1}$ and we are done.
\epr

\begin{center}
\textbf{Acknowledgements}
\end{center}

\noindent
We thank the Referee (of a journal) whose valuable and precise comments improved significantly the exposition and made precise several statements. In particular, the inclusion of the complete version of the proof of  Proposition 7, and the inclusion of Remark 1 in Section 3 were done following the Referee's suggestions.

\end{document}